\documentclass[draft]{article}
\usepackage{graphics}
\usepackage{amsmath}
\usepackage{amsfonts}
\usepackage{amssymb}
\usepackage{amsthm}

\usepackage{xtheorem}

\DeclareInstance{theoremstyle}{example}{std}
{
  pre-skip   = 1cm,
  post-skip  = 1cm,
  body-style = \ttfamily,
  head-style = \sffamily,
  post-head-action = \newline,
}

\newtheorem{plain}{Thm}{Theorem}[section]
\newtheorem{plain}{Prop}[Thm]{Proposition}
\newtheorem{plain}{Lem}[Thm]{Lemma}
{Corollary}
{Definition}
{Remark}
{Example}

\newcommand{\CC}{\mathbb C}
\newcommand{\HH}{\mathbb H}

\newcommand{\RR}{\mathbb R}
\newcommand{\ZZ}{\mathbb Z}

\begin{document}

\title{On harmonic and asymptotically harmonic homogeneous spaces}
\author{Jens Heber\thanks{Supported in part by DFG priority program ''Global Differential Geometry'' (SPP 1154)}
\\Mathematisches Seminar\\%
Christian-Albrechts-Universit\"at Kiel\\
Ludewig-Meyn-Str.~4\\D -- 24098 Kiel\\Germany}

\maketitle

\begin{abstract}
We classify noncompact homogeneous spaces which are Einstein 
and {\em asymptotically harmonic}. This completes the classification of Riemannian 
{\em harmonic spaces} in the homogeneous case: Any simply connected 
homogeneous harmonic space is flat, or rank-one symmetric, 
or a nonsymmetric Damek-Ricci space. Independently, Y.~Nikolayevsky 
has obtained the latter classification under the additional assumption 
of nonpositive sectional curvatures \cite{Ni2}.
\end{abstract}
\section{Introduction}
\label{intro}

A complete Riemannian manifold $(M,g)$ is called a {\em harmonic space} if about any point 
the geodesic spheres of sufficiently small radii are of constant mean curvature. 
If $M$ is noncompact and harmonic, then it is also {\em asymptotically harmonic}, 
that is, $M$ has no conjugate points and the mean curvature of 
its horospheres is constant (see Def.~\ref{def-harm}).

Clearly, any two-point homogeneous space is harmonic. In 1944, 
A.~Lichnerowicz asked the central question in this field: 
Is a harmonic space $M$ necessarily 
locally two-point homogeneous, that is, rank-one locally symmetric or flat?

Positive answers were given in 1944 by A.~Lichnerowicz for ${\rm dim}\,M=4$ \cite{Li}, 
in 1990 by Z.~Szab{\'o} for compact $M$ with 
finite fundamental group \cite{Sz}, 
in 1995 by G.~Besson, G.~Courtois, S.~Gallot for compact $M$ of negative sectional 
curvatures  \cite{BCG} (as a corollary of their proof of Gromov's minimal 
entropy conjecture, involving also results from 
\cite{BFL} and 
\cite{FL}), in 2002 by A.~Ranjan and H.~Shah 
for noncompact $M$ with minimal horospheres \cite{RS} (involving 
an observation of \cite{Ni1}), and 
in 2002 by Y.~Nikolayevsky \cite{Ni1} for ${\rm dim}\,M=5$.

However, in 1992, E.~Damek and F.~Ricci exhibited a class of noncompact 
harmonic spaces with sectional curvatures $K\le 0$ which are {\em homogeneous} but nonsymmetric 
\cite{DR}.

In this class, only the rank-one symmetric spaces 
have strictly negative curvature \cite{La}, \cite{Do2}. Any Damek-Ricci space 
admits a simply transitive, solvable Lie group of isometries $S$ 
whose commutator subgroup $[S,S]$ is $2$-step nilpotent and of codimension 
$1$ in $S$. Conversely, any {\em harmonic} homogeneous space modeled on 
such a Lie group $S$ is a Damek-Ricci space (see \cite{Dr2}, \cite{BPR}).

It is therefore natural to consider the following problems:

(a)\ Are simply connected harmonic spaces necessarily homogeneous?

(b)\ Classify all homogeneous harmonic spaces.

Question (a) appears to be widely open.

Recently, (b) was solved by Y.~Nikolayevsky \cite{Ni2} 
under the additional assumption of nonpositive sectional curvatures.

In this paper, we give an independent solution of (b) which does not 
require any a priori curvature assumption. In fact, we prove the following 
more general statement:

\begin{Thm}\label{class} Let $M$ be a noncompact, simply connected homo\-ge\-ne\-ous space. Then the 
following are equivalent:

$(i)$\,\ $M$ is {\em asymptotically harmonic} and Einstein.

$(ii)$\ $M$ is flat, or rank-one symmetric of noncompact type, or a nonsymmetric Damek-Ricci space.
\end{Thm}

We also obtain some algebraic restrictions without the assumption that 
$M$ is Einstein, i.~e.~of constant Ricci curvatures, see section~\ref{struct-proof}. 

Recall however, that every harmonic space is an Einstein manifold. Hence, Theorem~\ref{class} 
completes the classification of simply connected harmonic spaces in the homogeneous case:

\begin{Cor}\label{class-list} Let $M$ be a simply connected, {\em homogeneous harmonic space} of 
Ricci curvature $c\,.$ Then, up to scaling of the metric, $M$ is isometric to

$(a)$\ $S^n$, $\mathbb{C}P^k\,,$ $\mathbb{H}P^l\,,$ or $\mathbb{O}P^2\,,$ if $c>0.$

$(b)$\ $\mathbb{R}^n\,,$ if $c=0$.

$(c)$\ $\mathbb{R}H^n\,,\ \mathbb{C}H^k\,,\ \mathbb{H}H^l\,, \mathbb{O}H^2\,,$ or a 
nonsymmetric Damek-Ricci space, if $c<0.$ 
\end{Cor}
In fact, if $c>0$ holds, then it follows from the Bonnet-Myers Theorem and Z.~Szab\'o's 
Theorem \cite{Sz} that $M$ is rank-one symmetric of compact type. The case of $c=0$ 
is settled, since homogeneous Ricci flat spaces are necessarily flat \cite{AlK}. 
Any homogeneous Einstein space with $c<0$ is noncompact by a classical Bochner argument. 
If the space is harmonic in addition, then it has no conjugate points, is asymptotically 
harmonic (and Einstein) and Theorem \ref{class} applies directly.

Recalling from Corollary 1 of \cite{Heb3} (which is essentially based on \cite{AlC}) 
that nonsymmetric Damek-Ricci spaces 
do not admit any quotients of finite volume, we obtain

\begin{Cor}\label{fin-vol}
Let $M$ be a locally homogeneous, harmonic space of finite volume. Then, $M$ is 
rank-one locally symmetric or flat.
\end{Cor} 

Note that Riemannian products of nonflat harmonic spaces are not harmonic. They are, 
however, still D'Atri spaces, i.~e.~all of their (locally defined) geodesic symmetries 
preserve the volume element. Recall that the D'Atri property is in general quite flexible 
\cite{KPV}, for instance, any naturally reductive homogeneous space is a D'Atri space. 

However, in nonpositive curvature, we obtain the following rigidity result 
by combining the above results with Theorem~4.7 of \cite{Heb2}:

\begin{Cor}\label{DAtri} A homogeneous space $M$ with $K\le 0$ is a {\em D'Atri space}, iff 
it is isometric to a Riemannian product,
$$
\mathbb{R}^k\times T^l\times M_1\times \ldots \times M_m \times %
N_1\times \ldots \times N_n\,,
$$
where $T^l$ denotes a flat torus, each $M_i$ is irreducible symmetric of noncompact 
type and each $N_j$ is a nonsymmetric Damek-Ricci space. Any of the factors may be 
absent.
\end{Cor}

The organization of the paper is outlined in section~\ref{struct-proof}.

\section{Proof of Theorem~\ref{class}}\label{proof}
We recall and comment on the definition of harmonic and asymptotically harmonic 
spaces.

\begin{Def}\label{def-harm} Let $(M,g)$ denote a complete Riemannian manifold.

$(i)$ $M$ is called a {\em harmonic space}, if about any point 
the geodesic spheres of sufficiently small radii are of constant mean curvature.

$(ii)$ $M$ is called {\em asymptotically harmonic}, if 
$M$ has no conjugate points and the mean curvature of 
its horospheres is constant.
\end{Def}

For $M$ without conjugate points, {\em horospheres} in the 
universal covering ${\widetilde M}$ are defined as level sets of Busemann 
functions $b_v(x):=\lim_{t\to\infty} d(x,\gamma_v(t)) - t\,,\ x\in {\widetilde M}\,,\ %
v\in S{\widetilde M}.$ We call $b_v^{-1}(0)$ the {\em stable horosphere 
defined by} $v$. As for regularity of horospheres, see remark \ref{rem-harm} $(c)$ 
below.

For detailed references on harmonic spaces, we refer to \cite{Sz}. 
We recall some useful information.

\begin{Rem}\label{rem-harm}
$(a)$ Harmonicity is equivalent to the following (which explains the 
terminology): Any (locally defined) {\em harmonic function} satisfies the 
mean value property. It can be rephrased in terms of infinitely many conditions 
on the Riemann curvature tensor and its derivatives. These conditions are explicitly 
computable by a recursion formula and are named after A.~J.~Ledger (see 
e.~g.~\cite{Be}). For instance, harmonic spaces are Einstein, hence 
real analytic in normal coordinates, and condition $(i)$ extends to 
all radii except for those which correspond to conjugate points. 
It is well-known that the mean curvature constant of 
geodesic spheres in harmonic spaces depends only on the radius.

\smallskip

$(b)$ If $M$ and hence its Riemannian universal covering space 
${\widetilde M}$ is harmonic, then conjugate points occur at the 
same distance in every direction. Hence, either ${\widetilde M}$ is 
{\em compact}, or $M$ has {\em no conjugate points}. 
In the latter case, it follows from the Hadamard Cartan Theorem 
that ${\widetilde M}$ is diffeomorphic to a euclidean space and 
that any pair of points is joined by a unique (minimizing) geodesic.

\smallskip

$(c)$ In a general simply connected, complete Riemannian manifold $M$ without 
conjugate points, Busemann functions 
$b_v\,,\ v\in SM\,,$ 
are known to be of regularity $C^{1,1}$. We then define the ''mean curvature $m(v)$ of 
the stable horosphere $b_v^{-1}(0)$ at $\pi (v)$'' 
via the {\em stable Jacobi tensor} along $\gamma_v$:

Given $r>0$, consider the endomorphism fields 
$E_r(t) \in {\rm End}({\dot\gamma}_v(t)^{\perp})\,, t\in\RR\,,$ defined by 
the Jacobi equation $0=E_r''(t)+R(\cdot, {\dot\gamma}_v(t)){\dot\gamma}_v(t)\circ E_r(t)$ 
(''Jacobi tensors'') and the boundary conditions $E_r(0)={\rm id}, E_r(r)=0$. 
Note that $-E_r'(t)\circ E_r(t)^{-1}, t\in (-\infty ,r),$ is the second fundamental form of a 
sphere about $\gamma_v(r)$. The fields $t\mapsto E_r(t)$ converge locally uniformly on $\RR$ 
to the {\em stable Jacobi tensor} $t\mapsto E(t)\,,$ as $r\to\infty$. In fact, the second 
fundamental forms $U_r(v):=-E_r'(0)$ converge monotonically to $U(v):=-E'(0)$ (i.~e.~%
$U_R(v)-U_r(v)$ is negative definite for $r<R$), cf.~\cite{Gr}.

One then defines $m(v):={\rm trace}\,U(v)$. In particular, if $M$ has no conjugate points, 
then harmonicity implies asymptotic harmonicity.  

If $M$ is {\em asymptotically harmonic}, then $v\mapsto {\rm trace}\,U(v)$ is a 
nonnegative constant, and in particular continuous. By monotonical convergence and Dini's Theorem, 
the maps $v\mapsto {\rm trace}\,U_r(v)$ converge locally uniformly on $SM$, as $r\to\infty$. 
Hence, the tensor fields $U_r$ converge locally uniformly to $U$ \cite{Heb1}. Since $U_r(v)$ coincides 
with the Hessian of $b_{v,r}(x)=d(x,\gamma_v(r)) - r$ at $\pi (v)$, it follows 
that $b_{v,r}\,, r>0\,,$ converge in the local $C^2$-topology to $b_v$, 
as $r\to\infty$. Hence, in this case, Busemann functions are $C^2$ \cite{Esch}. 
Asymptotic harmonicity is equivalent to the existence of a constant $\alpha\ge 0$ 
such that $\Delta b ={\rm trace}_g{\rm Hess}\,b\equiv\alpha$ holds for any Busemann function 
$b:M\to\RR$. See \cite{HKS} for related discussion.

\smallskip

$(d)$ Suppose that $M$ is simply connected and asymptotically harmonic, 
say, ${\rm trace}\,U(v)\equiv\alpha\ge 0$, 
and that the $U_r$ converge {\em uniformly} to $U$ (e.~g.~if $M$ admits a compact quotient 
or if $M$ is homogeneous). Then, for any $\epsilon >0$, there exists an $R>0$ such 
that any sphere of radius at least $R$ has mean curvatures between $\alpha$ and 
$\alpha + \epsilon$. It follows that the logarithmic volume growth rate of $M$ 
equals $\alpha$, that is,
$$
\alpha = \lim_{r\to\infty} \frac{\log {\rm vol}\,B_r(p)}{r}\quad\mbox{for }p\in M.
$$
We say that $M$ has exponential (resp.~subexponential) volume growth, if $\alpha > 0$ 
(resp.~$\alpha =0$).

Note that the above limit exists for arbitrary complete $M$, 
provided that $M$ covers a compact manifold or that $M$ is homogeneous, 
as was proved by A.~Manning \cite{Ma}. If $M=(G/K, {\bar g})$ is an effective, Riemannian 
homogeneous space with connected Lie group $G$, we consider the Riemannian 
submersion $\pi :G\to G/K, h\mapsto [h]$ where $G$ is endowed with a suitable $G$-left invariant 
and $K$-right invariant metric $g$. Since $\pi^{-1}(B_r^{G/K}([e]))\subset %
B^G_{r+R}(e)\subset \pi^{-1}(B^{G/K}_{r+R}([e])$ holds for $R={\rm diam}(K)$, and 
${\rm vol}_g(\pi^{-1}(W))={\rm vol}_{\bar g}(W)\cdot {\rm vol}_g(K)$, we conclude 
that {\em subexponential growth of} $G/K$ can be characterized in the following 
equivalent ways: $(1)$ $G$ has subexponential growth, $(2)$ $G$ has 
polynomial growth, $(3)$ all ${\rm ad}_X\in {\rm End}(\frak{g}), X\in\frak{g},$ 
have only purely imaginary eigenvalues 
(cf.~sect.~6 of \cite{Pa} for proofs and the relevant references).
\end{Rem}

\subsection{Structure of the proof}\label{struct-proof}
We consider a simply connected homogeneous space $M$ which is {\em asymptotically 
harmonic}. The proof of Theorem~\ref{class} will be split into three parts. 
We refer to the corresponding sections for ideas of the proofs and further 
details.

\smallskip

$(i)$ In Proposition~\ref{prop-solv} we show that $M$ admits a simply 
transitive solvable Lie group of isometries $S$. Hence, $M$ is isometric to 
$S$, endowed with a suitable left invariant metric.

\smallskip

$(ii)$ In Proposition~\ref{prop-rank} we show 
that the commutator subgroup $N:= [S,S]$ has codimension $1$ in $S$, provided 
that $M$ has {\em exponential volume growth}. Hence, $S$ is a semidirect product 
of the normal nilpotent subgroup $N$ and a $1$-dimen\-sio\-nal complement 
$A\cong (\mathbb{R},+)$.

Note that in the case of {\em subexponential 
volume growth}, $S$ is a connected Lie group of subexponential and hence, 
polynomial volume growth (see remark \ref{rem-harm} $(d)$). 
All ${\rm ad}_X\in {\rm End}(\frak{s}), X\in\frak{s},$ 
have purely imaginary eigenvalues (cf.~also 
Cor.~\ref{expl-mean-curv}); in particular, $S$ is unimodular.\\
In the sequel, we will assume in addition that $M$ is an {\em Einstein} 
manifold. However by \cite{Do1}, any left invariant Einstein metric on a unimodular 
solvable Lie group is flat and then, so is $M$.

\smallskip

$(iii)$ In section~\ref{geom-char}, we recall facts about Damek-Ricci spaces. 
Finally, Theorem~\ref{asymp-harm-einst} completes the proof of 
Theorem~\ref{class}: We prove that any semidirect 
product $S=A\cdot N$ as above, endowed with an asymptotically harmonic 
Einstein metric, has constant curvature, or is a Damek-Ricci space.

\subsection{Transitive solvable isometry groups}
\label{sect-solv}

We prove the following (cf.~\cite{Wo} for related discussion):

\begin{Prop}\label{prop-solv} Every simply connected, homogeneous, 
{\em asymptotically harmonic} space $M$ 
admits a simply transitive solvable group $S$ of isometries.
\end{Prop}

\begin{proof} Let $G:={\rm Iso}_0(M)$ denote the identity component 
of the isometry group of $M$ and denote by $K:= G_p$ the isotropy 
subgroup of some point $p\in M$. 
Let ${\hat K}\supset K$ be a maximal compact subgroup of $G$. We prove 
that $K={\hat K}$:

Note that ${\hat K}$ is connected; moreover, $M_0:=G/{\hat K}$ is contractible, 
in fact, diffeomorphic to a euclidean space 
(cf.~Thm.~3.1 of \cite{I} or Thm.~3.1 of ch.~XV of 
\cite{Ho}). I follows that the canonical ${\hat K}$-principal fibre bundle 
$G\to G/{\hat K}= M_0$ is trivial (cf.~Cor.~10.3 of ch.~4 of \cite{Hu}), 
and that the associated projection fibre bundle $M=G/K\to G/{\hat K}=M_0$ 
with fibre ${\hat K}/K$ is trivial (see sect.~7 of ch.~4 of \cite{Hu}). 
In particular, $M$ is homotopy equivalent to the compact, connected manifold ${\hat K}/K$ 
and hence has nontrivial $\mathbb{Z}_2$-homology in dimension $n=\dim {\hat K}/K$.\\
On the other hand, $M$ is 
diffeomorphic to a euclidean 
space, since $M$ is simply connected with 
no conjugate points (Def.~\ref{def-harm}). It follows that ${\hat K}/K$ is $0$-dimensional. 
Hence, $K={\hat K}$ is maximal compact in $G$.

The assertion of the proposition follows from a standard Lie group 
argument, if the center $Z(G)$ is trivial \cite{Wo}, and requires a geometric 
argument in addition, if $Z(G)$ is nontrivial:

\smallskip

$(a)$ $Z(G)$ is trivial: Consider a Levi decomposition of the Lie 
algebra of $G$, say, $\frak{g}=\frak{l}\oplus\frak{r}=%
\frak{l}_c\oplus\frak{l}_{nc}\oplus\frak{r}$, where $\frak{r}$ denotes 
the radical of $\frak{g}$, and $\frak{l}$ denotes a semisimple complementary 
subgroup which we decompose into the direct sums $\frak{l}_c$ and $\frak{l}_{nc}$ 
of its compact ideals and its noncompact ideals, respectively. We denote the 
connected Lie subgroups, corresponding to this decomposition, by $L_c\,, L_{nc}$ 
and $R$, respectively. We also choose an Iwasawa decomposition 
$L_{nc}={\bar K}\cdot A\cdot N$ of the noncompact semisimple group and recall 
that $A$ is abelian and normalizes the nilpotent subgroup $N$, while 
${\bar K}$ contains the center ${\bar Z}:= Z(L_{nc})$ and ${\bar K}/{\bar Z}$ is compact.

Since $L_{nc}$ is a linear semisimple group, ${\bar Z}$ is finite and hence 
${\bar K}$ is compact. Since all maximal compact subgroups of $G$ are pairwise 
conjugate in $G$, we may assume that $K$ contains $L_c\cdot {\bar K}$. Hence, 
${\bar S}:= A\cdot N\cdot R$ is a solvable subgroup of $G$ which acts transitively 
on $M$. Since $M$ is simply connected, ${\bar S}$ contains a simply transitive 
(solvable) subgroup $S$, as asserted.

\smallskip

$(b)$ Suppose that $Z(G)$ is nontrivial; let ${\rm id}_M\not= \varphi\in Z(G)$: 
Then, the displacement function $p\mapsto d(p,\varphi (p))$ on $M$ is 
constant (''Clifford translation''), say, equal to $\alpha$ which is nonzero 
(since the isometry group acts effectively). 
For any $p\in M$, the geodesic $\gamma_{p\varphi (p)}$ through 
$p=\gamma_{p\varphi (p)}(0)$ and $\varphi (p)=\gamma_{p\varphi (p)}(\alpha)$ 
is $\varphi$-invariant: In fact for $0<t<\alpha$, we have $\alpha = %
d(\gamma_{p\varphi (p)}(t),\varphi\,\gamma_{p\varphi (p)}(t))%
\le d(\gamma_{p\varphi (p)}(t), \varphi (p))+d(\varphi (p), \varphi\,\gamma_{p\varphi (p)}(t))%
= \alpha$ and hence the broken geodesic from $\gamma_{p\varphi (p)}(t)$ via $\varphi (p)$ 
to $\varphi\,\gamma_{p\varphi (p)}(t)$ is a smooth geodesic. 

This provides a smooth foliation of $M$ by $\varphi$-invariant geodesics; 
choose one, say, $\gamma\,.$ For any $p\in M$, we obtain that
\begin{eqnarray*}
\lefteqn{b_{\gamma}(p)-b_{\gamma}(\varphi (p))=\lim_{t\to\infty} (d(p,\gamma(t))-d(\varphi (p),\gamma(t)))}\\
&=& \lim_{n\to\infty} (d(p,\varphi^{n+1}(\gamma (0)))-d(\varphi (p),\varphi^{n+1}(\gamma (0)))\\
&=& \lim_{n\to\infty} ((n+1) \alpha + b_{\gamma}(p)-n \alpha -b_{\gamma}(p))\\
&=& \alpha = d(p,\varphi (p))\,.
\end{eqnarray*}
Since Busemann functions are Lipschitz with Lipschitz constant $1$, 
it follows that the (unit length) gradient fields of $b_{\gamma}$ and 
$b_{\gamma^-}$ are tangent to the geodesic foliation, and hence, 
$b_{\gamma}+b_{\gamma^-}$ is constant. We conclude that 
$\Delta b_{\gamma}= -\Delta b_{\gamma^-}.$ However, since $M$ is 
asymptotically harmonic, it follows that $\Delta b_{\gamma} \equiv 0$ (cf.~remark \ref{rem-harm} $(c)$). 
Hence, all horospheres in $M$ are minimal, and $M$ has subexponential volume growth.

Now remark \ref{rem-harm} $(d)$ implies 
that $G$ has polynomial growth and hence, that all 
${\rm ad}_X\in {\rm End}(\frak{g}), X\in\frak{g},$ have only purely imaginary 
eigenvalues. In particular, $\frak{l}_{nc}$ is trivial. Finally, since 
$L_{nc}=\{{\rm id}_M\}$, it follows that the solvable radical $R$ 
acts transitively on $M$. As in step $(a)$, this completes the proof.
\end{proof}

\subsection{The algebraic rank of $S$}
\label{rank}

According to section~\ref{sect-solv}, any simply connected, homogeneous 
asymptotically harmonic space $M$ admits a simply transitive solvable group 
of isometries $S$. Hence, $M$ is isometric to $S$, endowed with a suitable 
left invariant metric $\langle\,,\rangle$ (''{\em solvmanifold}''). 
We now determine ${\rm codim}\,[S,S]$.  

\begin{Prop}\label{prop-rank} Let $M\cong (S,\langle\,,\rangle )$ be a solvmanifold which is 
asymptotically harmonic with {\em exponential volume growth}. Then, the commutator subgroup 
$N:= [S,S]$ has codimension $1$ in $S$.
\end{Prop}

We recall from $(ii)$ of section~\ref{struct-proof}, that asymptotically harmonic 
solv\-ma\-ni\-folds of {\em subexponential} volume growth are {\em flat}, provided that 
the metric is Einstein.

Proposition~\ref{prop-rank} yields directly the following

\begin{Cor}\label{cor-semidir}
$M\cong (A\cdot N,\langle\,,\rangle )$, semidirect 
product with $A\cong (\mathbb{R},+)$.
\end{Cor}

\begin{proof}[Proof of \ref{prop-rank}] We consider stable horospheres 
along geodesics perpendicular to $[S,S]\subset S$ and relate their shape operators 
to maximal solutions of certain algebraic Riccati equations (see (\ref{alg-ricc}) below). 
Results from control theory \cite{LR} allow one to compute the 
mean curvature of the horosphere in question (see (\ref{mean-curv-a})). Since mean 
curvatures are constant by assumption, the claim follows.

\smallskip

Consider the metric Lie algebra $(\frak{s}, \langle\,,\rangle )$ 
of left invariant vector fields on $S$ 
and its orthogonal vector space decomposition $\mathfrak{s}=\mathfrak{a}\oplus\mathfrak{n}$ 
(note that $\mathfrak{a}$ need not be a subalgebra!). Choose a unit vector $A\in\mathfrak{a}$. 
If we decompose ${\rm ad}_A=D_A+S_A$ into its symmetric and skew symmetric part 
(w.~r.~t.~$\langle\,,\rangle$), then Jacobi operator and covariant derivative 
along $A$ are given by $R_A:=R(\cdot ,A)A=-D_A^2- [D_A,S_A]$ and $\nabla_A=S_A$ on $\frak{s}$\,, 
respectively (cf.~e.~g.~\cite{AW1}). Note that $\nabla_AA=0$.

We compute the mean curvature of stable horospheres $H(t)$ along $\gamma_A(t)$. 
Their shape operators $L(t)\in {\rm End}({\dot\gamma}_A(t)^{\perp})$ satisfy the Riccati equation
\begin{eqnarray}\label{ricc}
L'(t)+L^2(t)+R(t)=0\quad (L(t)=L(t)^T)
\end{eqnarray}
(where $R(t)=R(\cdot\,,{\dot\gamma}_A(t)){\dot\gamma}_A(t)$). 
Since $\gamma_A$ is a one-parameter subgroup of $S$, 
it follows that, in an orthonormal basis of left invariant vector fields, $L(t)$ is constant, say, 
$L(t)\equiv L_0=L_0^T\in {\rm End}(A^{\perp})$. More precisely, writing $Y(t):=Y\circ\gamma_A$ 
for any $Y\in\frak{s}$, we obtain
\begin{eqnarray*}
\begin{aligned}
0&=L'(t)Y(t)+L(t)^2Y(t)+R(t)Y(t)\\
&=(\nabla_A (L_0Y))(t)-(L_0 \nabla_AY)(t) +(L_0^2Y)(t)+(R_AY)(t)\\
&=([S_A,L_0] Y+L_0^2Y-D_A^2Y-[D_A,S_A] Y)(t)\,.
\end{aligned}
\end{eqnarray*}
If we write $L_0=-D_A-X$, then we obtain that $X=X^T$ is a 
solution of the algebraic (matrix) Riccati equation
\begin{eqnarray}\label{alg-ricc}
X^2+X\circ{\rm ad}_A+{\rm ad}_A^T\circ X=0\,.
\end{eqnarray}
Since any solution of (\ref{alg-ricc}) corresponds to a solution of (\ref{ricc}) 
defined globally on $\RR$ and $L$ is the minimal among all such solutions 
(see e.~g.~Proposition 3' of \cite{EO}), 
it follows that $X$ is the unique {\em maximal solution} of (\ref{alg-ricc}). 
We recall from \cite{LR} that for this solution, $-{\rm ad}_A-X$ has eigenvalues 
of nonpositive real parts. Now since (\ref{alg-ricc}) is equivalent to
\begin{eqnarray*}
\begin{aligned}
\left(\begin{array}{cc}I&0\\X&I\end{array}\right)^{-1}%
\left(\begin{array}{cc}-{\rm ad}_A&-I\\0&{\rm ad}_A^T\end{array}\right)%
\left(\begin{array}{cc}I&0\\X&I\end{array}\right)&\\[0.5ex]
= \left(\begin{array}{cc}-{\rm ad}_A-X&-I\\0&{\rm ad}_A^T+X\end{array}\right)&\,,%
\end{aligned}
\end{eqnarray*}
it follows that any eigenvalue of $-{\rm ad}_A-X$ is also an eigenvalue of $-{\rm ad}_A$ or of 
$+{\rm ad}_A^T\,.$ We conclude that
\begin{eqnarray}\label{mean-curv-a}
{\rm trace}\,L_0={\rm trace}\,(-{\rm ad}_A-X)=%
-\sum\limits_{\sigma}\vert {\rm Re}(\sigma )\vert\,,
\end{eqnarray}
where the sum is taken over all eigenvalues $\sigma$ of ${\rm ad}_A$, with multiplicities.

Since $\frak{s}$ is solvable, it follows from Lie's Theorem 
that the eigenvalues of ${\rm ad}_X$ depend linearly on $X\in\frak{s}$. 
On the other hand, since $M$ is asymptotically harmonic of exponential volume growth, 
the expression in (\ref{mean-curv-a}) is nonzero and independent of $A\in\frak{a},\Vert A\Vert =1$. 
This is impossible unless $\mathfrak{a}$ is one-dimensional.
\end{proof}

The proof includes the following explicit information
\begin{Cor}\label{expl-mean-curv}
Let $M\cong (S,\langle\,,\rangle )$ be a solvmanifold and let 
$\gamma_A\,, A\in\frak{s}, \Vert A\Vert =1\,,$ 
denote a geodesic perpendicular to $[S,S]$ at $e$. If $\gamma_A$ contains 
no conjugate points, then the trace of the stable Riccati solution 
along $\gamma_A$ is equal to
$$
-\sum\limits_{\sigma}\vert {\rm Re}(\sigma )\vert\,,
$$
where the sum is taken over all eigenvalues $\sigma$ of ${\rm ad}_A$, with multiplicities.
\end{Cor}

\subsection{Geometric characterization of Damek-Ricci spaces}
\label{geom-char}

In this section, we consider a simply connected homogeneous space $M$ which is 
asymptotically harmonic and Einstein. We prove in Theorem~\ref{asymp-harm-einst} 
that $M$ is isometric to a Damek-Ricci space (Def.~\ref{def-heis-dr}), 
unless it is a space of constant sectional curvature $K\equiv c\le 0.$

Following Proposition~\ref{prop-solv}, $M$ is isometric to a solvable Lie group $S$, 
endowed with a left invariant metric $\langle\,,\rangle$. By $(ii)$ of 
section~\ref{struct-proof}, either $M$ is flat, or $M$ has exponential 
volume growth. If the latter condition is satisfied, then $N=[S,S]$ has codimension 
$1$ in $S$, as follows from Proposition~\ref{prop-rank}; we will restrict our 
attention to this algebraic structure.

\medskip

A Lie algebra $\frak{n}$ is called $2$-step nilpotent, it its derived algebra 
$[\frak{n},\frak{n}]$ lies in the center $\frak{z}:=\frak{z}(\frak{n})$. If 
$\frak{n}$ is $2$-step nilpotent and is endowed with a positive definite scalar 
product $\langle\,,\rangle$, then we consider the orthogonal complement 
$\frak{v}$ of $\frak{z}$ in $\frak{n}$. Following A.~Kaplan \cite{Ka}, the Lie 
bracket of $\frak{n}$ defines an endomorphism
$$
j:\frak{z}\to \frak{so}(\frak{v},\langle\,,\rangle )\,,\quad 
\langle j(Z) V, W\rangle :=\langle [V,W], Z\rangle\quad\mbox{for } Z\in\frak{z}\,, V,W\in\frak{v}\,,
$$
with values in the skew symmetric endomorphisms of $(\frak{v},\langle\,,\rangle)$. 
Note that, conversely, any such endomorphism defines a $2$-step nilpotent Lie 
bracket on $\frak{n}=\frak{v}\oplus\frak{z}$ (with center containing $\frak{z}$).

\begin{Def}\label{def-heis-dr}
$(i)$ A metric $2$-step nilpotent Lie algebra $(\mathfrak{n},[,],\langle\,,\rangle )$ is 
said to be {\em of Heisenberg type} (or {\em of type}~H) if its $j$-map satisfies
$$
j(Z)^2=- {\rm id}_{\mathfrak{v}}\quad\mbox{ for all }Z\in\mathfrak{z}, \Vert Z\Vert =1.
$$

$(ii)$ A solvable Lie group with left invariant metric $(S,\langle\,,\rangle )$ 
is called a {\em Damek-Ricci space} if $\frak{n}:=[\frak{s},\frak{s}]$ is $2$-step 
nilpotent of Heisenberg type, and $\frak{n}^{\perp}$ 
is spanned by a unit vector $H$, such that ${\rm ad}_H$ equals $1/2\cdot {\rm id}$ on $\frak{v}$ 
and $1\cdot {\rm id}$ on $\frak{z}$.
\end{Def}

\begin{Rem}\label{rem-heis-dr}
Condition $(i)$ implies that $j$ extends to an algebra homomorphism between 
the Clifford algebra $Cl(\frak{z},\langle\,,\rangle )$ and ${\rm End}(\frak{v})$, 
i.~e.~Heisenberg type algebras with $l$-dimensional center are in $1$-$1$-correspondence 
with modules over $Cl_l$ (in the notation of \cite{LM}). Among all Heisenberg 
type algebras, the Iwasawa subalgebras of the real semisimple Lie algebras of 
split rank one are of particular 
interest: The corresponding Damek-Ricci spaces are the rank-one symmetric 
spaces of noncompact type (and nonconstant curvature), that is, the hyperbolic 
spaces over $\CC\,, \HH$ and the hyperbolic Cayley plane.
\end{Rem}

Damek-Ricci spaces satisfy $K\le 0$ and are Einstein \cite{Bo}, 
in fact, harmonic \cite{DR}. In this class, only the rank-one symmetric spaces 
have strictly negative curvature \cite{La}, \cite{Do2}. For a detailed discussion 
of Damek-Ricci spaces, we refer to \cite{BTV} (for homogeneous spaces 
with $K\le 0\,,$ cf.~\cite{Hei}, \cite{Al}, \cite{AW1}, \cite{AW2}; 
for noncompact homogeneous Einstein spaces, cf.~\cite{Heb3}).

Rigidity 
results for Damek-Ricci spaces have been obtained under suitable 
algebraic assumptions: Suppose that $S$ is a solvable Lie group endowed 
with a left invariant {\em harmonic} Riemannian metric $g$. It was proved 
in \cite{BPR} and in \cite{Dr2} (under further weak algebraic restrictions) 
that $(S, g)$ is isometric to a Damek-Ricci space 
provided that the commutator subgroup $[S,S]$ is two-step nilpotent and 
of codimension one.

Our proof of \ref{asymp-harm-einst} yields that the above-mentioned 
additional a priori algebraic assumptions on $[S,S]$ are not necessary. Moreover, 
asymptotic harmonicity and the Einstein condition are sufficient as the 
geometric requirements. 
The proof starts with explicit calculations of certain geodesics and Jacobi 
operators by M.~Druetta \cite{Dr1} and was in part inspired by those.

\begin{Thm}\label{asymp-harm-einst} 
Let $S$ denote a simply connected solvable Lie group with 
commutator subgroup $N=[S,S]$ of codimension $1$. If $S$ 
admits a nonflat, asymptotically harmonic Einstein metric $\langle\,,\rangle$, 
then up to scaling of the metric, $(S,\langle\,,\rangle)$ is constantly 
curved or isometric to a Damek-Ricci space.
\end{Thm}

\begin{proof}
We consider certain geodesics $\gamma_Z$, emanating from $e$, tangent to the center of $N$ 
($(a), (b)$). In $(b)$, we determine all stable Jacobi fields along $\gamma_Z$ 
explicitly in terms of hypergeometric functions. We then obtain formula (\ref{mean-curv}) 
for the mean curvature of stable horospheres along $\gamma_Z$ (cf.~$(c)$). The formula involves 
a product of expressions in hypergeometric functions (\ref{def-h}) which, by asymptotic harmonicity, 
has to be constant. An analytic continuation argument (cf.~$(d)$) then shows that, in fact, every 
factor of this product is constant. This yields strong restrictions on the algebraic 
structure of $S$ (cf.~$(e)$), as required.

\smallskip

$(a)$ We decompose the Lie algebra of $S$ orthogonally as 
$\mathfrak{s}=\langle H\rangle \oplus \mathfrak{v}\oplus\mathfrak{z}$ 
where $\frak{n}=[\frak{s}, \frak{s}]=\frak{v}\oplus\frak{z}\,, \frak{z}=\frak{z}(\frak{n})\,,$ 
and $H$ is a unit vector 
(note that $\mathfrak{n}=\mathfrak{v}\oplus\mathfrak{z}$ is nilpotent, 
since $\frak{s}$ is solvable, but not a priori $2$-step).

Since ${\rm codim}\,N=1$, it follows that $(S,\langle\,,\rangle)$ is a 
noncompact homogeneous Einstein space of {\em standard type}, as defined 
in \cite{Heb3}. According to the structure theory developed in \cite{Heb3}, 
the subspace $\frak{v}=\frak{n}\cap\frak{z}^{\perp}$ is ${\rm ad}_H$-invariant; moreover 
(possibly after passing to an isometric metric Lie algebra, see Theorem~4.10 
of \cite{Heb3}),
$$
{\rm ad}_H:\mathfrak{v}\to \mathfrak{v},\quad {\rm ad}_H:\mathfrak{z}\to \mathfrak{z}
$$
are both self-adjoint with {\em positive} eigenvalues
$$
\rho_1\le \ldots\le\rho_k\ \,\mbox{resp.}~\mu_1\le \ldots\le\mu_l=:\lambda
$$
(possibly after replacing $H$ by $-H$). We will not make use of the fact (cf.~\cite{Heb3}) 
that the eigenvalues are integers, multiplied by a constant factor.

The eigenspace of ${\rm ad}_H$ corresponding 
to the maximal eigenvalue lies in the center $\frak{z}$ of $\frak{n}$. Hence, the maximal 
eigenvalue equals $\lambda$, we have $\rho_k<\lambda\,$ and, after suitably rescaling the metric, we 
may and will assume that $\lambda =1$.

We choose an arbitrary unit length eigenvector $Z\in \mathfrak{z}\,, %
{\rm ad}_HZ=\lambda Z\,.$ As in the beginning of section~\ref{geom-char}, 
we define a skew symmetric endomorphism $j(Z)\in\mathfrak{so}(\mathfrak{v})$ by 
$\langle [V,{\tilde V}], Z\rangle_{\mathfrak{z}}=\langle j(Z)V,{\tilde V}\rangle_{\mathfrak{v}}$ 
(note that the collection of all $j(Z)\,, Z\in\frak{z}\,,$ gives only partial information 
about the Lie bracket, unless $\frak{n}$ is $2$-step nilpotent).

For later reference, we note that
\begin{eqnarray}\label{j-ad-commute}
[j(Z)^2,{\rm ad}_H]=0\quad\mbox{on }\frak{v}\,.
\end{eqnarray}
In fact, if ${\rm ad}_HV=\rho V\,, {\rm ad}_HW=\rho' W\,, V,W\in\frak{v}\,,$ 
then $[V,W]$ lies in the $(\rho+\rho')$-eigenspace of ${\rm ad}_H$. Hence, 
$\langle j(Z)V,W\rangle =\langle [V,W], Z\rangle$ vanishes unless 
$\rho +\rho' =1$. It follows that $j(Z)$ maps any $\rho$-eigenspace of ${\rm ad}_H\vert_{\frak{v}}$ 
to its $(1-\rho)$-eigenspace, and that $j(Z)^2$ leaves any eigenspace of 
${\rm ad}_H$ on $\frak{v}$ invariant.

\smallskip

In the sequel, we will prove that $\mu_j=1$ holds for all $j$, $\rho_i=1/2$ holds for all $i$, 
and that $j(Z)^2=- {\rm id}_{\mathfrak{v}}\,.$ Since this holds for all $Z\in\frak{z}\,, %
\Vert Z\Vert =1$, we find that either $\frak{v}=0\,, {\rm ad}_H\vert_{\frak{z}}={\rm id}_{\frak{z}}$ and 
$(S,\langle\,,\rangle)$ is isometric to real hyperbolic space, or that $\frak{n}$ is nonabelian, 
i.~e.~$\frak{v}\not= 0\,,$ and that $(S,\langle\,,\rangle)$ is isometric 
to a Damek-Ricci space. This observation will complete the proof.

\medskip

$(b)$ We recall from Lemma 1 of \cite{Dr1} that the geodesic with initial velocity $Z$ 
is given by $\gamma_Z(t)=\exp (\tanh (t) Z)\cdot\exp (-\log\cosh (t) H)$ (in terms of 
the Lie group exponential map of $S$). 
In terms of the left invariant vector fields $H$ and $Z$ on $S$, the velocity vector 
field is given by ${\dot\gamma}_Z(t)=(-\tanh (t) H+\frac{1}{\cosh (t)} Z)\vert_{\gamma_Z(t)}\,.$ 

We determine {\em stable} (in fact, all) {\em Jacobi fields} along $\gamma_Z$ 
explicitly in terms of {\em hypergeometric functions} 
(with parameters given by ${\rm ad}_H$ and $j(Z)$) and in a suitable basis of 
left invariant vector fields.

\smallskip

Among all Jacobi 
fields, we identify the stable ones as those whose norm is {\em bounded for} $t\in[0,\infty )$. 
In fact, since $R(\cdot ,H)H=-{\rm ad}^2_H$ is negative definite on $H^{\perp}$, 
we conclude that $R(t):=R(\cdot,{\dot\gamma}_Z(t)){\dot\gamma}_Z(t)$ is 
negative definite on ${\dot\gamma}_Z(t)^{\perp}$ for $t$ large enough, say, $t\ge t_0$. 
Hence, $t\mapsto\Vert J(t)\Vert$ is convex on $[t_0,\infty )$ for any Jacobi field. 
Since $(S,\langle ,\rangle )$ has no conjugate points, any $v\in Z^{\perp}$ defines unique 
Jacobi fields $J_r$ along $\gamma_Z$ such that $J_r(0)=v\,, J_r(r)=0$ and the stable field 
$J=\lim_{r\to\infty} J_r$. Obviously, $t\mapsto \Vert J(t)\Vert$ is bounded 
(in fact, decreasing on $[t_0, \infty )$), and $J$ is the only Jacobi field 
with initial value equal to $v$ which is bounded on $[0,\infty ).$

We show that the vector space of stable Jacobi fields 
is the direct sum of three subspaces, containing $(b1):$ stable Jacobi fields 
which are everywhere tangent to $\langle H, Z\rangle$, $(b2):$ those which are 
everywhere tangent to $(\frak{z}\cap Z^{\perp})$, and $(b3), (b4):$ those tangent to 
$\frak{v}$.

We calculate covariant derivatives of left invariant vector fields, 
using the Koszul formula $2 \langle \nabla_XY,W\rangle =%
-\langle X,[Y,W]\rangle -\langle Y,[X,W]\rangle + \langle W,[X,Y]\rangle\,, X,Y,W\in\frak{s}$ (see 
e.~g.~section 1.3 of \cite{Dr1}): For instance, since ${\rm ad}_H$ is self-adjoint, 
$\nabla_HY=0$ holds for all $Y\in\frak{s}$. Using that 
${\dot \gamma}_Z(t)=(-\tanh (t) H+\frac{1}{\cosh (t)} Z)\vert_{\gamma_Z(t)}$, a routine 
calculation yields the following identities which hold for all left invariant 
vector fields $Z^*\in\frak{z}$ with $Z^*\perp Z$ and for all $V\in\frak{v}$:
\begin{eqnarray}\label{nabla-Z*-V}
\nabla_{{\dot \gamma}_Z(t)}Z^*=0\,,\quad 
\nabla_{{\dot \gamma}_Z(t)}V=\frac{-1}{2 \cosh (t)} (j(Z)V)\vert_{\gamma_Z(t)}\,.
\end{eqnarray} 

$(b1)$ Since $\nabla_HZ=\nabla_HH=0$ and $\nabla_ZH=-Z\,, \nabla_ZZ=H$, it follows 
that $H$ and $Z$ span a totally geodesic $2$-dimensional subalgebra 
of constant curvature $-1$. We obtain a $1$-dimensional space of orthogonal 
stable Jacobi fields along $\gamma_Z$ which are tangent to $\exp (\langle H, Z\rangle )$.

\smallskip

$(b2)$ Suppose that ${\rm dim}\,\frak{z}\ge 2\,.$ Let $Z^*\in\frak{z}\,, Z^*\perp Z\,,$ 
be a unit length eigenvector of ${\rm ad}_H$, 
say, ${\rm ad}_HZ^*=\mu Z^*\,.$ Recall that $0 < \mu \le \lambda =1$ 
and that the left invariant vector field $Z^*$ is 
parallel along $\gamma_Z$ (cf.~(\ref{nabla-Z*-V})).

Note that  
$R(t)=\frac{-1}{\cosh^2 (t)}({\rm ad}_H+\sinh^2 (t){\rm ad}_H^2)$ on 
$(\frak{z}\cap Z^{\perp})\oplus {\rm ker}(j(Z))$ (cf. Prop.~1 of \cite{Dr1}). In particular, $R(t)Z^*=%
\frac{-\mu-\sinh^2 (t)\mu^2}{\cosh^2 (t)}Z^*$. We obtain a $2$-dimensional 
vector space of Jacobi fields of the form $t\mapsto h(t)Z^*\vert_{\gamma_Z(t)}$ 
where $h$ is any solution of
$$
0=h''(t)+\frac{-\mu-\sinh^2 (t)\mu^2}{\cosh^2 (t)} h(t)
$$
or, equivalently, $h(t)= e^{-\mu t}g(t)$ and
$$
0=g''(t)-2 \mu g'(t) +\frac{\mu^2-\mu}{\cosh^2 (t)} g(t)\,.
$$
Substituting $g(t)=u(z(t))$ with $z=z(t)=\frac{1-\tanh (t)}{2}\in (0,1)$, we obtain
\begin{eqnarray}\label{hypgeom-cent}
0=z (1-z) u''(z) + (1+\mu - 2 z) u'(z) - \mu (1 - \mu) u(z)\,,
\end{eqnarray}
the hypergeometric equation 
$z (1-z) f''(z)+(c-(a+b+1) z) f'(z) - a b f(z) = 0$ 
with parameters $a= \mu\,, b=1-\mu\,, c= 1+\mu\,.$ The hypergeometric function
\begin{eqnarray}\label{hyp-geom-fct}
F(a,b;c;z)=1+\sum_{k=1}^{\infty}\frac{\prod_{l=0}^{k-1}(a+l) \prod_{l=0}^{k-1}(b+l)}{\prod_{l=0}^{k-1}(c+l)}%
\cdot\frac{z^k}{k!}
\end{eqnarray}
provides a solution of (\ref{hypgeom-cent}) which is regular at $z=0$ (cf.~sect.~2.1.1 
of \cite{EMOT}). 
This yields a Jacobi field 
$t\mapsto e^{-\mu t}F(\mu,1-\mu; 1+\mu;\frac{1-\tanh (t)}{2}) Z^*$ along $\gamma_Z$ 
whose norm is bounded for $t\in [0,\infty)$, that is, a stable Jacobi field.

Note that a 
linearly independent solution of (\ref{hypgeom-cent}) is given by $u(z)=z^{-\mu}$. This yields a 
multiple of the Jacobi vector field $t\mapsto\cosh^{\mu}(t) Z^*$, restriction of 
the right invariant {\em Killing vector field} defined by $Z^*\in T_eS$.

\smallskip

$(b3)$ According to (\ref{j-ad-commute}), the kernel of $j(Z)\in\frak{so}(\frak{v})$ 
is ${\rm ad}_H$-invariant. Suppose that ${\rm ker}(j(Z))$ is nontrivial. 
Choose a unit length eigenvector $V^*\in {\rm ker}(j(Z))$ of ${\rm ad}_H$, 
say, ${\rm ad}_HV^*=\rho^* V^*\,.$ Recall that $0 < \rho^* < \lambda =1$ 
and that $V^*$ is parallel along $\gamma_Z$ (cf.~(\ref{nabla-Z*-V})). Moreover, $R(t)V^*=%
\frac{-\rho^*-\sinh^2 (t){\rho^*}^2}{\cosh^2 (t)}V^*$, and arguing as in $(b2)$, 
we obtain a stable Jacobi field 
$t\mapsto e^{-\rho^* t}F(\rho^*,1-\rho^*; 1+\rho^*;\frac{1-\tanh (t)}{2}) V^*$ along $\gamma_Z$.

\smallskip

$(b4)$ The orthogonal complement of ${\rm ker}(j(Z))$ in $\frak{v}$ decomposes 
into an orthogonal direct sum of two-dimensional subspaces which are ${\rm ad}_H$- and 
$j(Z)$-invariant. Any such subspace admits an orthonormal basis $\{ V,{\tilde V}\}$ 
such that the matrix representations are
\begin{eqnarray}\label{matrix-j-ad}
{\rm ad}_H=\left( \begin{array}{cc}\rho & 0\\0 & 1-\rho\end{array}\right)\,,\quad %
j(Z)=\left( \begin{array}{cc}0 & -\theta\\+\theta & 0\end{array}\right)\,,\quad %
0<\rho\le 1/2\,,\quad\theta > 0\,.
\end{eqnarray}
Let $V(t):=V(\gamma_Z(t))$ and ${\tilde V}(t):={\tilde V}(\gamma_Z(t))$. 
We show that stable Jacobi fields 
along $\gamma_Z$ with initial value in $\langle V, {\tilde V}\rangle$ are of 
the form $t\mapsto f(t) V(t)+g(t) {\tilde V}(t)$ with functions $f, g$ given 
explicitly in terms of hypergeometric functions (\ref{stable-fields}):

\smallskip

Note that $\frac{DV}{dt}(t)=\frac{-\theta}{2 \cosh (t)}{\tilde V}(t)$ and 
$\frac{D{\tilde V}}{dt}(t)=\frac{+\theta}{2 \cosh (t)}V(t)$ 
(cf.~(\ref{nabla-Z*-V}) and (\ref{matrix-j-ad})). 
In particular, the left invariant vector fields $V$ and ${\tilde V}$ span a 
parallel bundle along $\gamma_Z$. For later reference, we conclude that 
$\frac{D^2V}{dt^2}(t)=\frac{+\theta \sinh (t)}{2 \cosh^2 (t)}{\tilde V}(t)-%
\frac{\theta^2}{4 \cosh^2 (t)}V(t)$ and 
$\frac{D^2{\tilde V}}{dt^2}(t)=\frac{-\theta \sinh (t)}{2 \cosh^2 (t)}V(t)-%
\frac{\theta^2}{4 \cosh^2 (t)}{\tilde V}(t).$ 
On the other hand, the Jacobi operator acts on left invariant 
vector fields $\in\frak{v}$ as 
$R(t)=\frac{1}{\cosh^2 (t)} \{ -\frac{1}{4} j(Z)^2-{\rm ad}_H-\sinh^2 (t) {\rm ad}_H^2 %
-\sinh (t) j(Z)\circ (1/2\cdot {\rm id}-{\rm ad}_H)\}$ (Prop.~1 of \cite{Dr1}), 
hence leaving $\langle V,{\tilde V}\rangle$ invariant. Therefore, we obtain a 
$4$-dimensional space of Jacobi fields along $\gamma_Z$ of the form 
$t\mapsto f(t) V(t)+g(t) {\tilde V}(t)$ which solve
\begin{eqnarray}\label{jac-eqn}
\begin{aligned}
0&=\left(\frac{D^2}{dt^2}+R(t)\right)(f(t) V(t)+g(t) {\tilde V}(t))\\
=&\left( f''(t)+\frac{\theta}{\cosh (t)}g'(t)-\frac{\rho+\sinh^2(t)\rho^2}{\cosh^2(t)}f(t)%
+\frac{(\theta\rho-\theta) \sinh (t)}{\cosh^2 (t)}g(t)\right) V(t)\\
+&\left( g''(t)+\frac{-\theta}{\cosh (t)}f'(t)-\frac{(1-\rho)+\sinh^2(t)(1-\rho)^2}{\cosh^2(t)}g(t)%
+\frac{\theta\rho\sinh (t)}{\cosh^2 (t)}f(t)\right) {\tilde V}(t).
\end{aligned}
\end{eqnarray}
In order to solve (\ref{jac-eqn}), consider
$$
A(t)=\tanh (t)\left( \begin{array}{cc}\rho & 0\\0 & 1-\rho\end{array}\right)\,,\quad %
B(t)=A(t)+\frac{1}{\cosh (t)}\left( \begin{array}{cc}0&-\theta\\{+\theta} & 0\end{array}\right)\,.
$$
Then (\ref{jac-eqn}) is equivalent to either of the following two equations
\begin{eqnarray}\label{jac-eqn-2}
\begin{aligned}
0&=\left(\frac{d}{dt}+A(t)\right)\circ \left(\frac{d}{dt}-B(t)\right)%
\left(\begin{array}{c}f\\g\end{array}\right)(t)\,,\\
0&=\left(\frac{d}{dt}+B^T(t)\right)\circ \left(\frac{d}{dt}-A(t)\right)%
\left(\begin{array}{c}f\\g\end{array}\right)(t)\,.
\end{aligned}
\end{eqnarray}
The vector space of solutions $\left({f\atop g}\right)(t)$ of (\ref{jac-eqn}) 
hence contains the $2$-dimensional subspaces 
${\rm ker}\,(\frac{d}{dt}-A(t))$ and ${\rm ker}\,(\frac{d}{dt}-B(t))$ 
and is in fact spanned by these (since any element $\left({f\atop g}\right)(t)$ 
in the intersection satisfies 
$(B(t)-A(t))\left({f\atop g}\right)(t)=0$ and hence vanishes).

As in $(b2)$, $(b3)$, it is convenient to write down these solutions 
in terms of $z=z(t)=\frac{1-\tanh(t)}{2}\,.$ 
First, 
${\rm ker}\,(\frac{d}{dt}-A(t))$ is spanned by 
$t\mapsto \left({\cosh^{\rho}(t)\atop 0}\right)=\left({(4z(1-z))^{-\rho /2}\atop 0}\right)$ 
and 
$t\mapsto \left( {0\atop \cosh^{1-\rho}(t)}\right)=%
\left( {0\atop (4z(1-z))^{(\rho -1)/2}}\right)$ (these solutions correspond to 
right invariant {\em Killing vector fields} restricted to $\gamma_Z$).
It remains to exhibit ${\rm ker}\,(\frac{d}{dt}-B(t))$:\\

To that end, let $u:(0,1)\to\mathbb{R}$ be any solution of the 
hypergeometric equation
\begin{eqnarray}\label{hypgeom-u}
0=z (1-z) u''(z) + (\rho - 2 \rho z) u'(z) + \theta^2 u(z)\,.
\end{eqnarray}
(Only in step $(d)$ of the proof, we will consider holomorphic extensions.) 
As usual, we consider the coefficients $a,b,c$ of (\ref{hypgeom-u}) which, in terms of 
$\rho$ and $\theta$, are given by $c=\rho\,, a+b+1 =2\rho\,, a b =-\theta^2$. 
Equivalently,
\begin{eqnarray}\label{hypgeom-u-2}
0=\cosh^{-2}(t) u''(z(t))+4\rho\tanh (t) u'(z(t))%
+4\theta^2 u(z(t))\,.
\end{eqnarray}
Using (\ref{hypgeom-u-2}), a straightforward calculation yields that
\begin{eqnarray}\label{ker-D-B}
\left(\begin{array}{c}f\\g\end{array}\right)(t)=%
\left(\begin{array}{c}-\cosh^{-\rho}(t) u'(z(t))\\%
2\theta \cosh^{1-\rho}(t) u(z(t))\end{array}\right)=
\left(\begin{array}{c}-(4 z(1-z))^{\rho /2} u'(z)\\%
2\theta (4z(1-z))^{(\rho -1)/2} u(z)\end{array}\right)
\end{eqnarray}
solves $\left({f\atop g}\right) '(t)=B(t)\left({f\atop g}\right)(t)$ and hence, 
(\ref{jac-eqn-2}) resp.~(\ref{jac-eqn}) (where we again abbreviate $z=z(t)=\frac{1-\tanh (t)}{2}$).

As solutions of (\ref{hypgeom-u}), we use (cf.~2.8 (20), (22), 2.9 (1), (2) of \cite{EMOT})
\begin{eqnarray}\label{examples-h}
\begin{aligned}
u_1(z)&=F(a,b;c;z)\\
u'_1(z)&=\frac{a b}{c} F(a+1,b+1;c+1;z)\\[1ex]
u_2(z)&=z^{1-c} F(1+a-c,1+b-c;2-c;z)\\
u'_2(z)&=(1-c) z^{-c} F(1+a-c,1+b-c;1-c;z)\\
&= (1-c) (z (1-z))^{-c} F(-a,-b;1-c;z)\,.
 \end{aligned}
\end{eqnarray}
The solutions, given by (\ref{ker-D-B}), $u=u_1$ and $u=u_2$, respectively, 
are not bounded for $t\in [0,\infty)$, but adding suitable elements from 
${\rm ker}\,(\frac{d}{dt}-A(t))$ we obtain bounded solutions:

We consider the following special solutions of (\ref{jac-eqn}) resp.~(\ref{jac-eqn-2}).
\begin{eqnarray}\label{stable-fields}
\begin{aligned}
\left(\begin{array}{c}f\\g\end{array}\right) (t)&=B_{\rho}(z)\cdot C_{\rho ,\theta}(z)\cdot %
\left(\begin{array}{c}v_1\\v_2\end{array}\right)\mbox{ for }z=z(t)=\frac{1-\tanh (t)}{2}\,,\ %
v_i\in\RR\,,\\%
B_{\rho}(z)&=\left(\begin{array}{cc}(4 z(1-z))^{-\rho /2} & 0\\0 & (4 z(1-z))^{(\rho -1) /2}\end{array}%
\right) ,\\%
C_{\rho , \theta}(z) &= %
\left(\begin{array}{cc}-(4 z(1-z))^{\rho}u'_1(z) & -(4 z(1-z))^{\rho}u'_2(z)+4^{\rho} (1-\rho)\\%
2\theta u_1(z)-2\theta & 2\theta u_2(z)\end{array}\right) . %
\end{aligned}
\end{eqnarray}
Note that $B_{\rho}(z)\cdot C_{\rho ,\theta}(z)\to 0$ as $z=z(t)\to 0$ resp.~$t\to\infty$. 
This follows since $0<\rho\le 1/2$ and from the definition of $u_i(z)$ (cf.~(\ref{examples-h})) 
in terms of hypergeometric functions. One uses the fact that hypergeometric functions are of the form 
$F(\ldots ;z)=1+z G(\ldots ;z)$ with $G(\ldots ;z)$ analytic at $z=0$ (cf.~(\ref{hyp-geom-fct})). 
Hence, the solutions in (\ref{stable-fields}) are bounded 
for $t\in [0,\infty)$ and, by the argument in the beginning 
of step $(b)$, they provide all {\em stable Jacobi fields} 
with initial value in $\langle V, {\tilde V}\rangle$.

For later reference, consider the Wronskian $W(z):=W(u_1,u_2)(z)$ of (\ref{hypgeom-u}), 
that is,  
$W(z)=u_1(z) u'_2(z)-u_2(z) u'_1(z)\Rightarrow W'(z)=\frac{2\rho z-\rho}{z(1-z)} W(z) %
\Rightarrow W(z)=const\cdot (z(1-z))^{-\rho}\,;$ plugging in (\ref{examples-h}) and 
letting $z\to 0$, we obtain $const =1-\rho\,,$ 
and hence,
\begin{eqnarray}\label{det-E-2}
\lefteqn{{\rm det}\,B_{\rho}(z)\cdot C_{\rho ,\theta}(z)}\nonumber\\
&=&-4^{\rho } (1-\rho) \theta\cdot\sqrt{\frac{z}{1-z}}\cdot\left( %
\frac{u_1(z)+F(-a,-b;1-c;z)-2}{z}\right)\nonumber\\
&=&const_{\rho ,\theta}\cdot e^{-t}\cdot \left( \frac{F(a,b;\rho;z)+F(-a,-b;1-\rho ;z)-2}{z}\right) .
\end{eqnarray}
\medskip

$(c)$ Recall the Busemann function $b_Z(x)=\lim_{t\to\infty} d(x,\gamma_Z(t))-t$ on 
$(S,\langle\,,\rangle)$ and the stable horosphere $b_Z^{-1}(-t)$ through $\gamma_Z(t)$. 
We compute its mean curvature $m(t)$ at $\gamma_Z(t)$ explicitly in terms of 
hypergeometric functions (\ref{mean-curv}). By asymptotic harmonicity, 
it then follows that a certain product, whose factors involve hypergeometric functions 
(\ref{def-h}) is constant.

\smallskip

To that end, consider the stable Jacobi tensor $t\mapsto E(t)\in {\rm End}({\dot\gamma}_Z(t)^{\perp})$ 
as defined in remark \ref{rem-harm} $(c)$. For any parallel field $t\mapsto v(t)$, 
the stable Jacobi field with initial value $v(0)$ is given by $t\mapsto E(t) v(t)$. 

Choose an orthonormal basis of $Z^{\perp}\in T_eS$ which consists of eigenvectors of ${\rm ad}_H$ as 
considered in $(b1)-(b4)$, say, $H$ and bases $\{ Z^*_j\}_j$ of $\frak{z}\cap Z^{\perp}$, 
$\{V^*_k\}_k$ of $\frak{v}\cap {\rm ker}(j(Z))$, and $\bigcup_i \{ V_i, {\tilde V}_i\}$ for 
$\frak{v}\cap {\rm ker}(j(Z))^{\perp}$; that is,
\begin{eqnarray}\label{jac-data}
\begin{aligned}
{\rm ad}_HZ^*_j&=\mu_j Z^*_j\quad &{\rm ad}_HV^*_k&=\rho^*_k V^*_k\\
{\rm ad}_HV_i&=\rho_i V_i \quad &{\rm ad}_H{\tilde V}_i&=(1-\rho_i) {\tilde V}_i\\
j(Z)V_i&=+\theta_i {\tilde V}_i\quad &j(Z){\tilde V}_i&=-\theta_i V_i
\end{aligned}
\end{eqnarray}
where $0<\mu_j \le 1\,,\ 0<\rho^*_k<1\,,\ 0<\rho_i\le 1/2\,,\ \theta_i>0$ for all $i, j, k\,.$ 
We consider the matrix representation of $E(t)$ w.~r.~t.~the parallel orthonormal frame field 
along $\gamma_Z$ corresponding to this basis. We obtain a block 
matrix structure with $(1\times 1)$-blocks with entries $e^{-t}$ (cf.~$(b1)$), 
$e^{-\mu_j t}F(\mu_j,1-\mu_j; 1+\mu_j;z)$ (cf.~$(b2)$), 
$e^{-\rho^*_k t}F(\rho^*_k,1-\rho^*_k; 1+\rho^*_k;z)$ (cf.~$(b3)$) and $(2\times 2)$-blocks 
$E_i(t)$ for each $i$ (cf.~$(b4)$). The entries of $E_i(t)$ 
are coefficients 
of two stable Jacobi fields, written in a {\em parallel} frame field, while 
$B_{\rho_i}(z)\cdot C_{\rho_i ,\theta_i}(z)$ as defined in (\ref{stable-fields}) 
encodes two different stable fields 
in terms of the {\em left invariant} fields $V_i$ and ${\tilde V}_i$. Hence, 
the matrices are transformed into each other by multiplication 
with a suitable constant invertible matrix from the right 
and a (time dependent) orthogonal matrix from the left. We conclude that 
${\rm det}\,E_i(t)=const\cdot {\rm det}\,(B_{\rho_i}(z(t))\cdot C_{\rho_i ,\theta_i}(z(t)))$ 
holds for each $i$. 
It then follows from (\ref{det-E-2}) that
\begin{eqnarray}\label{det-E}
{\rm det}\,E(t)= const\cdot e^{-\nu t}\cdot h\left(\frac{1-\tanh (t)}{2}\right)
\end{eqnarray}
where $\nu =1+\sum_j\,\mu_j +\sum_k\,\rho^*_k+\sum_i\,1 ={\rm trace}\,{\rm ad}_H\,,$ 
and the function $h$ is defined by
\begin{eqnarray}\label{def-h}
\begin{aligned}
h(z):= &\prod_j\,F(\mu_j,1-\mu_j; 1+\mu_j;z)\cdot\prod_k\,F(\rho^*_k,1-\rho^*_k; 1+\rho^*_k;z)\\
 &\cdot\prod_i\,\frac{F(a_i,b_i;\rho_i;z)+F(-a_i,-b_i;1-\rho_i;z)-2}{z}\,.
\end{aligned}
\end{eqnarray}
Here, the $a_i, b_i$ are given by $a_i+b_i+1 =2\rho_i\in (0,1]\,, a_i b_i =-\theta_i^2<0\,.$ 
We may and will assume that $a_i < 0 < b_i$ and hence, that $b_i\le \vert a_i\vert < b_i+1\,.$

The mean curvature $m(t)$ of the horosphere $b_Z^{-1}(-t)$ at $\gamma_Z(+t)$ is then 
given by
\begin{eqnarray}\label{mean-curv}
m(t)=-\frac{d}{dt}\,\log\,\vert {\rm det}\,E(t)\vert =%
{\rm trace}\,{\rm ad}_H - \frac{d}{dt}\,\log\,\vert h\left(\frac{1-\tanh (t)}{2}\right)\vert\,.
\end{eqnarray}
The second summand on the right hand side of (\ref{mean-curv}) tends to $0$, as $t\to\infty$, 
since $h(0)=\prod_i\,(\frac{a_i b_i}{\rho_i}+\frac{a_i b_i}{1-\rho_i})\not= 0$. 
Since $S$ endowed with the left invariant metric $\langle\,,\rangle$ is, 
by assumption, asymptotically harmonic, it follows that $t\mapsto m(t)$ is 
constant equal to ${\rm trace}\,{\rm ad}_H\,.$ Hence, the function $h$ is 
constant on $(0,1)$ (the image of $t\mapsto \frac{1-\tanh (t)}{2}$). 
We conclude that $z\mapsto h(z)$ is constant.

\medskip

$(d)$ We prove that each of the factors of $h(z)$ is constant:

Recall that every hypergeometric function (and hence, every factor of $h(z)$) 
can be continued analytically along any path which avoids $0\,, 1$ and $\infty$. 
Analytic continuation along a closed path yields another solution of the 
corresponding hypergeometric equation close to the endpoint.

For each factor $f(z)$ of $h(z)$, we investigate its analytic continuation ${\tilde f}(z)$ 
along a simple positive loop $\beta$, based at $1/2$, around $1$. We will prove that either 
$\lim_{t\to 0} \vert {\tilde f}(t)\vert =\infty$, or $f(z)={\tilde f}(z)$ is a polynomial. 
Since the product of all $f$'s is constant (and different from $0$), it then follows 
that all of the factors are constant functions. In step $(e)$, we will derive 
restrictions on the defining parameters $\mu_j\,, \rho^*_k\,, \rho_i$ and $\theta_i$.

\smallskip

Given $a,b,c\in\RR$ with $c\notin\ZZ$, consider the hypergeometric equation
$$
0=z(1-z)\,u''(z)+(c-(a+b+1) z)\,u'(z) - a b\, u(z)
$$
and its solutions $u_1(z)=F(a,b;c;z)$ and 
$u_2(z)=z^{1-c} F(1+a-c,1+b-c;2-c;z)$ (with the principal branch of $z^{1-c}$) 
which are linearly independent and uniquely defined, say, on 
$\CC\backslash ((-\infty ,0]\cup [1,\infty))$.  
Recall from section~2.7.1 of \cite{EMOT}, that 
the analytic continuation of $u_1$ along $\beta$ equals
\begin{eqnarray}\label{holonomy}
\begin{aligned}
{\tilde u_1} &= B_{11} u_1 + B_{12} u_2\,,\quad\mbox{where}\\
B_{11}&=1-2 i e^{i \pi (c-a-b)}\frac{\sin (\pi a) \sin (\pi b)}{\sin (\pi c)}\,,\\
B_{12}&=- 2 i \pi e^{i \pi (c-a-b)}%
\frac{\Gamma (c) \Gamma (c-1)}{\Gamma (c-a) \Gamma (c-b)\Gamma (b)\Gamma (a)}\,.
\end{aligned}
\end{eqnarray}
Note that $B_{12}= 0\,,$ iff at least one of the numbers $a, b, c-a, c-b$ lies in 
$\ZZ_0^{-}$.

\smallskip

$(d1)$ If $a=\mu_j\,, b=1-\mu_j\,, c=1+\mu_j$ and $0<\mu_j< 1$, then $B_{12}\not= 0\,.$ 
Moreover, $u_2(z)=z^{-\mu_j}$ and hence, $\lim_{z\to 0}\vert {\tilde u_1}(z)\vert = \infty$ 
(where ${\tilde u_1}(z)$ denotes the analytic continuation of $u_1(z)$ along $\beta$). 

If $\mu_j =1\,,$ then $1\equiv u_1(z) = {\tilde u_1}(z)\,.$

\smallskip

$(d2)$ If $a=\rho^*_k\,, b=1-\rho^*_k\,, c=1+\rho^*_k$ and $0<\rho^*_k< 1$, then 
$\lim_{z\to 0}\vert {\tilde u_1}(z)\vert = \infty$, compare $(d1)$.
\smallskip

$(d3)$ Factors of $h(z)$ of the third type are of the form
$$
u(z)=\frac{F(a,b;c;z)+F(-a,-b;1-c;z)-2}{z}\,,
$$
where $0<c\le 1/2\,,\ a<0<b$ and $a+b+1=2 c$. In order to exhibit 
the analytic continuation ${\tilde u}(z)$ of $u(z)$ along $\beta$, 
we consider all three summands of $u(z)$ separately, and find that
$$
{\tilde u}(z)=A \frac{1}{z}+ B \frac{1}{z^c} + C \frac{1}{z^{1-c}} + g(z)
$$
where $A,B,C\in\mathbb{C}$ are explicitly computable in terms of (\ref{holonomy}) 
and $g$ is analytic, say, in $\CC\backslash ((-\infty ,0]\cup [1,\infty))$ 
and bounded near $0$. Hence, 
$\lim_{z\to 0}\vert {\tilde u}(z)\vert = \infty\,,$ unless $A=0$ and 
either $B=C=0$ or $c=1/2\,, B=-C\,.$ The latter conditions imply restrictions 
on $a,b,c$ as follows:

If $c=1/2$, then $a+b=2 c-1=0$; it follows that
\begin{eqnarray*}
\begin{aligned}
A&=-2 i e^{i \pi (c-a-b)}\frac{\sin (\pi a) \sin (\pi b)}{\sin (\pi c)}%
-2 i e^{i \pi (1-c+a+b)}\frac{\sin (\pi a) \sin (\pi b)}{\sin (\pi (1-c))}\\%
&=-2 i e^{i \pi (1-c)}\frac{\sin (\pi a) \sin (\pi b)}{\sin (\pi c)}%
-2 i e^{i \pi c}\frac{\sin (\pi a) \sin (\pi b)}{\sin (\pi (1-c))}
\end{aligned}
\end{eqnarray*}
equals $-4 \sin^2 (\pi b)$ which vanishes iff $b\in\ZZ^+$; hence, $u(z)=2 \frac{F(-b,b;1/2;z)-1}{z}$ 
is a polynomial of degree $b-1$ (cf.~(\ref{hyp-geom-fct})).

If $c\not= 1/2\,,$ then $A=0$ implies that $a$ or $b$ are integers. 
Since $c$ is not an integer, it follows that $c-b=1-c+a$ is not an integer. 
But then, using the $B_{12}$-vanishing criterion in (\ref{holonomy}), 
$B=0$ implies that $a$ is an integer, 
while $C=0$ yields that $-b$ is an integer. Finally, 
since $0<a+b+1=2 c\le 1\,,$ it follows that $a=-b\,,\ c=1/2\,,$ a contradiction. 

\medskip

$(e)$ As proved in $(d)$, asymptotic harmonicity implies that every factor 
of $h(z)$ is constant. We conclude that factors involving $\rho^*_k$ are not 
present (cf.~$(d2)$), that all $\mu_j$ are equal to $1$ (cf.~$(d1)$), 
that 
$a_i=-1\,, b_i=+1\,, \rho_i=1/2$ holds for all $i$ (cf.~$(d3)$) and 
hence, $\theta_i=\sqrt{-a_i b_i}=+1\,.$

We can reformulate this as follows: $\frak{v}\cap {\rm ker}(j(Z))$ is 
trivial, ${\rm ad}_H$ equals $1\cdot {\rm id}$ on $\frak{z}$ and 
$1/2\cdot {\rm id}$ on $\frak{v}$, and $j(Z)^2=-{\rm id}_{\frak{v}}$. 
As explained in the end of step $(a)$, this completes the proof.
\end{proof}

{\bf Acknowledgement.} The present paper is part of a joint project 
with H.~Shah and G.~Knieper. We are indebted to DFG for the support 
within the priority program 
''Global Differential Geometry'' (SPP 1154). I am grateful to 
H.~Shah, G.~Knieper and U.~Abresch for helpful discussions and 
would like to express my gratitude for the hospitality at the 
Mathematical Institute at the Ruhr-Universit\"at Bochum. I would also like 
to thank the referee whose detailed comments were very helpful improving the 
readability of the paper.

\end{document}